\documentclass[graybox]{svmult}

\RequirePackage{type1cm}        \RequirePackage{makeidx}         \RequirePackage{graphicx}        \RequirePackage{multicol}        \RequirePackage[bottom]{footmisc}

\RequirePackage{amssymb}
\RequirePackage{xfrac}

\RequirePackage{newtxtext}       \usepackage{enumitem}

\RequirePackage{calrsfs}

\RequirePackage{amsfonts}

\RequirePackage[most]{tcolorbox}

\RequirePackage{bbm}

\RequirePackage{graphicx}

\RequirePackage{booktabs}

\RequirePackage{url}

\RequirePackage[hidelinks]{hyperref}

\RequirePackage{color}

\RequirePackage{wasysym}
\RequirePackage{listings}
\RequirePackage{caption}
\RequirePackage{subcaption}

\RequirePackage{stmaryrd}

\RequirePackage{lineno}

\RequirePackage{xparse}

\RequirePackage{float}

\RequirePackage{multirow}

\RequirePackage{cellspace}

\RequirePackage{shellesc}

\RequirePackage[ruled,vlined,linesnumbered]{algorithm2e}

\SetCommentSty{xCommentSty}

\RequirePackage{tikz}
\usetikzlibrary{calc, automata, positioning, arrows, external, matrix, shapes.misc}
\tikzset{cross/.style={cross out,thick,draw=black,minimum size=2*(#1-\pgflinewidth), inner sep=0pt, outer sep=0pt},
cross/.default={2.5pt}}

\RequirePackage{pgfplots}
\pgfplotsset{compat=1.16}
\usepgfplotslibrary{fillbetween}
\usepgfplotslibrary{groupplots}
\usepgfplotslibrary{external}
\DeclareMathAlphabet{\pazocal}{OMS}{zplm}{m}{n}
\DeclareMathAlphabet\bpazocal{OMS}{cmsy}{b}{n}

\providecommand{\vertiii}[1]{{\left\vert\kern-0.15ex\left\vert\kern-0.15ex\left\vert #1
    \right\vert\kern-0.15ex\right\vert\kern-0.15ex\right\vert}}

\NewDocumentCommand{\curlii}{sO{}m}
{
	\IfBooleanTF{#1}
    {\dgalext{#3}}
    {\dgalx[#2]{#3}}
}

\NewDocumentCommand{\dgalext}{m}{  \sbox0{    \mathsurround=0pt     $\left\{\vphantom{#1}\right.\kern-\nulldelimiterspace$  }  \sbox2{\{}  \ifdim\ht0=\ht2
    \{\kern-.625\wd2 \{#1\}\kern-.625\wd2 \}  \else
    \left\{\kern-.7\wd0\left\{#1\right\}\kern-.7\wd0\right\}  \fi
}

\NewDocumentCommand{\dgalx}{om}{  \sbox0{\mathsurround=0pt$#1\{$}  \sbox2{\{}  \ifdim\ht0=\ht2
    \{\kern-.625\wd2 \{#2\}\kern-.625\wd2 \}  \else
    \mathopen{#1\{\kern-.7\wd0 #1\{}
    #2
    \mathclose{#1\}\kern-.7\wd0 #1\}}
  \fi
}

\tikzset{
  partial ellipse/.style args={#1:#2:#3}{
    insert path={+ (#1:#3) arc (#1:#2:#3)}
  }
}

\providecommand{\to}{\widetilde{o}}

\definecolor{blackmy}{RGB}{38, 70, 83}
\definecolor{bluemy}{RGB}{39, 125, 161}
\definecolor{greenmy}{RGB}{42, 167, 143}
\definecolor{yellowmy}{RGB}{233, 196, 106}
\definecolor{brownmy}{RGB}{244, 162, 97}
\definecolor{redmy}{RGB}{249, 65, 68}
\usepackage{color, colortbl}

\usepackage{lineno}

\makeatletter
\newcommand{\thickhline}{\noalign {\ifnum 0=`}\fi \hrule height 1pt
    \futurelet \reserved@a \@xhline
}
\newcolumntype{"}{@{\hskip\tabcolsep\vrule width 1pt\hskip\tabcolsep}}
\def\blfootnote{\gdef\@thefnmark{}\@footnotetext}
\makeatother

\makeindex

\begin{document}

\title*{Integrating Additive Multigrid with Multipreconditioned Conjugate Gradient Method}
\author{Hardik Kothari\orcidID{0000-0003-0110-5384} and
	\\Maria Giuseppina Chiara Nestola\orcidID{0000-0002-5700-0306} and
	\\Marco Favino\orcidID{0000-0002-1253-9164} and
	\\Rolf Krause\orcidID{0000-0001-5408-5271}}
\institute{Hardik Kothari \at Euler Institute, Universit\`{a} della Svizzera italiana, Switzerland, \email{hardik.kothari@usi.ch}
\and
Maria Giuseppina Chiara Nestola \at Euler Institute, Universit\`{a} della Svizzera italiana, Switzerland, \email{nestom@usi.ch}
\and
Marco Favino\at Euler Institute, Universit\`{a} della Svizzera italiana, Switzerland, \email{marco.favino@usi.ch}
\at UniDistance Suisse, \email{marco.favino@fernuni.ch}
\and
Rolf Krause
\at King Abdullah University of Science and Technology, \email{rolf.krause@kaust.edu.sa}
\at UniDistance Suisse, \email{rolf.krause@fernuni.ch}
}
\authorrunning{Hardik Kothari et al.}
\maketitle
\abstract*{Due to its optimal complexity, the multigrid (MG) method is one of the most popular approaches for solving large-scale linear systems arising from the discretization of partial differential equations.
	However, the parallel implementation of standard MG methods, which are inherently multiplicative, suffers from increasing communication complexity.
	In such cases, the additive variants of MG methods provide a good alternative due to their inherently parallel nature, although they exhibit slower convergence.
	This work combines the additive multigrid method with the multipreconditioned conjugate gradient (MPCG) method.
	In the proposed approach, the MPCG method employs the corrections from the different levels of the MG hierarchy as separate preconditioned search directions.
	In this approach, the MPCG method updates the current iterate by using the linear combination of the preconditioned search directions, where the optimal coefficients for the linear combination are computed by exploiting the energy norm minimization of the CG method.
	The idea behind our approach is to combine the $A$-conjugacy of the search directions of the MPCG method and the quasi $H_1$-orthogonality of the corrections from the MG hierarchy.
	In the numerical section, we study the performance of the proposed method compared to the standard additive and multiplicative MG methods used as preconditioners for the CG method.}

\abstract{Due to its optimal complexity, the multigrid (MG) method is one of the most popular approaches for solving large-scale linear systems arising from the discretization of partial differential equations.
	However, the parallel implementation of standard MG methods, which are inherently multiplicative, suffers from increasing communication complexity.
	In such cases, the additive variants of MG methods provide a good alternative due to their inherently parallel nature, although they exhibit slower convergence.
	This work combines the additive multigrid method with the multipreconditioned conjugate gradient (MPCG) method.
	In the proposed approach, the MPCG method employs the corrections from the different levels of the MG hierarchy as separate preconditioned search directions.
	In this approach, the MPCG method updates the current iterate by using the linear combination of the preconditioned search directions, where the optimal coefficients for the linear combination are computed by exploiting the energy norm minimization of the CG method.
	The idea behind our approach is to combine the $A$-conjugacy of the search directions of the MPCG method and the quasi $H_1$-orthogonality of the corrections from the MG hierarchy.
	In the numerical section, we study the performance of the proposed method compared to the standard additive and multiplicative MG methods used as preconditioners for the CG method.}

\section{Problem definition}
In this work, we focus on modeling a boundary value problem fundamental to the mathematical modeling of diffusion processes.
We consider a Lipschitz domain $\Omega \subset \mathbb{R}^d$, $d=2$, with the boundary $\partial \Omega$.
The boundary is decomposed into the Dirichlet boundary $\Gamma_D$ and the Neumann boundary $\Gamma_N = \partial\Omega \setminus \Gamma_D$.
The problem is defined as: given a symmetric diffusion tensor $\boldsymbol{K}:\Omega \to \mathbb{R}^{d \times d}_{\mathrm{sym}}$ and $f:\Omega \to \mathbb{R}$, find $u: \Omega \to \mathbb{R}$ such that
\begin{equation}
	\begin{aligned}
		-\nabla \cdot \boldsymbol{K}(x) \nabla u & = f   &  & \quad \text{in} \ \Omega,   \\
		u                                        & = u_D &  & \quad \text{on} \ \Gamma_D, \\
		\nabla u \cdot \boldsymbol{n}            & = g_N &  & \quad \text{on} \ \Gamma_N.
	\end{aligned}
	\label{eq:diffusion}
\end{equation}
We consider specifically anisotropic diffusion, which is characterized by the diffusion tensor, $\boldsymbol{K}(x) = \begin{pmatrix} k_{xx}(x) & 0 \\ 0 & k_{yy}(x) \end{pmatrix}$, where {${k_{xx}(x) \ll \ k_{yy}(x)}$}.

In the context of groundwater flow, \eqref{eq:diffusion} models the fluid movement through subsurface layers.
We define these layers as a set of equidimensional fractures, denoted as $\Omega_f$, while the matrix domain is given as $\Omega_m = \Omega \setminus \Omega_f$.
The flow in fractured media is of particular interest due to its relevance in geophysics, specifically for oil and gas recovery, geothermal energy, and pollution control.
Here, the diffusion tensor within the matrix region is defined as $\boldsymbol{K}(x) = k_{m} \mathbbm{1}_d$, for all $x \in \Omega_m$, and within the fracture region as $\boldsymbol{K}(x) = k_{f} \mathbbm{1}_d$, for all $x \in \Omega_f$, where $\mathbbm{1}_d$ represents the identity matrix of dimension $d$.

To numerically solve~\eqref{eq:diffusion}, the finite element method (FEM) is utilized for discretization.
This approach effectively transforms the continuous diffusion problem into an algebraic system of equations, typically formulated as $\boldsymbol{A} \boldsymbol{x} = \boldsymbol{b}$.

\section{Solution method}
This section will discuss our approach to utilizing the additive multilevel method within the multipreconditioned conjugate gradient method.
Before discussing our approach, we briefly introduce the additive multigrid method and then discuss our solution strategy.

\subsection{Additive multigrid method}
The multigrid method we usually employ in practice is a multiplicative multigrid method.
As in the $V$/$W$-cycle, the global correction is computed using the successive application of the smoothing steps on each level while going first from the fine-to-coarse levels and then from the coarse-to-fine levels~\cite{7_hackbusch_multi-grid_1986}. 
In contrast, the additive MG method allows for the application of the smoothing operations in parallel, completely independent of each other~\cite{2_bastian_additive_1998}.
Thus, the additive MG method presents a benefit over traditional multiplicative MG techniques due to its inherently parallel nature, which can be advantageous in modern high-performance computing environments~\cite{5_chen_multigrid_2001, 6_darwish_parallelization_2008}.

To realize a multigrid method, we assume that the multilevel hierarchy of FE spaces is constructed either in a geometric manner, i.e., by constructing a hierarchy of nested meshes, or in an algebraic manner, i.e., by employing some agglomeration methods.
We denote this hierarchy of nested FE spaces as ${ \pazocal{V}_0 \subset \pazocal{V}_1 \subset \cdots \subset \pazocal{V}_L }$, where $L$ denotes the finest level and $0$ denotes the coarsest level of the multilevel hierarchy.
On this hierarchy of FE spaces, we construct the prolongation operators that can transfer the functions from a coarse level to the finest level, given as
\(
P_\ell^L: \pazocal{V}_\ell \to \pazocal{V}_L, \text{for } \ell \in \{0,1, \ldots, L\}.
\)
Similarly, we also define the restriction operators that transfer the residuals from the finest level to a coarse level, given as
\(
R_L^\ell: \pazocal{V}^\ast_L \to \pazocal{V}^\ast_\ell, \text{for } \ell \in \{0,1, \ldots, L\},
\)
where for the linear MG method, the restriction operators can be constructed as the transpose of the prolongation operator, e.g., $R_L^\ell = (P_\ell^L)^\top$, and $\pazocal{V}^\ast_\ell$ denotes the dual space of $\pazocal{V}_\ell$.
Also, for $\ell = L$, the transfer operators degenerate to the identity operator.
We note that these transfer operators for the additive MG method can be constructed explicitly or computed as composite operators by successively applying the transfer operators defined between two consecutive levels. 
Additionally, we define smoothing operators on each level of the multilevel hierarchy except the coarsest level as:
\(
\pazocal{S}_\ell: \pazocal{V}^\ast_\ell \to  \pazocal{V}_\ell, \text{for } \ell \in \{1,2, \ldots, L\},
\)
where the smoothing iterations for given initial guess $\boldsymbol{c}_0$ on level $\ell$ are given as
\begin{equation}
	\pazocal{S}_\ell^{\nu}(R_L^\ell \boldsymbol{b}) := (\pazocal{I} - \pazocal{M}_\ell^{-1} \pazocal{A}_\ell )^\nu \boldsymbol{c}_0 + \sum_{i=0}^{\nu-1}(\pazocal{I} - \pazocal{M}_\ell^{-1} \pazocal{A}_\ell )^i \pazocal{M}_\ell^{-1}(R_L^\ell \boldsymbol{b}),
\end{equation}
where $\pazocal{M}_\ell$ and  $\pazocal{I} - \pazocal{M}_\ell^{-1} \pazocal{A}_\ell$ denote a preconditioner and an iteation matrix induced by the smoother.
On the coarsest level, we define an approximate inverse of the coarse level operator $\pazocal{A}_0$ as $\pazocal{S}_0: \pazocal{V}_0^\ast \to \pazocal{V}_0$.
As a single step of the smoothing iteration may not be sufficient to smooth the error associated with the high-frequency components of the error, it becomes essential to perform multiple smoothing iterations, which we define as $\nu$.
With these operators, we can define the update of the additive multigrid method as
\begin{equation}
	\boldsymbol{x}_{k+1} = \boldsymbol{x}_{k} + P_0^L \pazocal{S}_0 (R^0_L \boldsymbol{r}_k)+ \sum_{\ell=1}^{L} P_\ell^L \pazocal{S}^\nu_\ell (R^\ell_L \boldsymbol{r}_k),
	\label{eq:additiveMG}
\end{equation}
where $\boldsymbol{r}_k = \boldsymbol{b} - \boldsymbol{A} \boldsymbol{x}_k$.
As we can observe, in \eqref{eq:additiveMG} the corrections from all levels are equally weighted, and the method might also diverge if the corrections are not appropriately damped.
Due to this reason, the additive MG method is usually employed as a preconditioner, where the optimal damping parameter is computed by the outer iterative method, usually a Krylov subspace method.

In practice, the additive MG method has cheaper synchronization costs than the standard multiplicative MG method on parallel architectures~\cite{1_bastian_load_1998}.
However, the convergence of the additive MG-PCG method remains slower than its multiplicative counterpart, which prohibits the usage of the additive MG-PCG method as a solution method for large-scale problems.
We aim to improve the convergence of the overall methodology by combining the additive MG method and the multipreconditioned CG method.

\subsection{Multipreconditioned Conjugate Gradient method (MPCG)}

The multipreconditioned CG method was introduced as a variant of a CG method with the ability to employ multiple preconditioners~\cite{4_bridson_multipreconditioned_2006}.
In the MPCG approach, multiple preconditioning operators are employed at each iteration to construct more search directions, and rather than summing up the search directions with equal weights, the coefficients for the linear combination for search directions are computed using the energy minimizing property.
By employing the multiple preconditioners, the MPCG method constructs the iterates in the generalized Krylov space while maintaining the $A$-conjugacy and orthogonality property of the standard CG method.
However, using multiple preconditioners eliminates the three-term recursion feature of the CG method, necessitating the adoption of the truncated CG method.
This approach requires the storage of search directions from the previous $m$ iterations.
These stored directions facilitate the generation of a brief recurrence relation to construct new search directions, ensuring their orthogonality to the preceding $m$ directions.

In this work, we integrate the MPCG method with the additive MG approach.
A key aspect of our methodology is the utilization of corrections from various levels of the multigrid hierarchy, each serving as a distinct preconditioner.
For the iteration $k$, the additive MG preconditioned residuals are given as
\[
	\boldsymbol{P}_k = [P^L_0 \pazocal{S}_0 (R_L^0 \boldsymbol{r}_k) \ \vert \ P^L_1 \pazocal{S}_1^\nu (R_L^1 \boldsymbol{r}_k) \ \vert   \cdots \  \vert \ P^L_L \pazocal{S}^{\nu}_L (R_L^L \boldsymbol{r}_k)] \in \mathbb{R}^{n\times(L+1)}.
\]
Using these search directions, the step sizes for each direction in the MPCG method is constructed using the following formulation,
\[
	\mathbb{R}^{L+1} \ni \boldsymbol{\alpha} = (\boldsymbol{P}_{k}^\top \boldsymbol{A} \boldsymbol{P}_{k})^{-1}(\boldsymbol{P}_{k}^\top \boldsymbol{r}_{k}),
\]
which are employed to construct the global search direction, given as
\begin{equation}
	\boldsymbol{P}_k \boldsymbol{\alpha} = \boldsymbol{p}_{k} = \alpha_0 P_0^L  \pazocal{S}_0 (R_L^0 \boldsymbol{r}_k ) + \sum_{\ell=1}^{L}\alpha_\ell P_\ell^L  \pazocal{S}^\nu_\ell (R_L^\ell \boldsymbol{r}_k).
	\label{eq:linear_combination}
\end{equation}

With this approach, we can repurpose the smoothing operations on the different levels of the multilevel hierarchy as actions of different preconditioners.
This approach leverages different subspaces to reduce the error in different frequency of the spectrum, thereby constructing a more robust preconditioning approach.
In the context of the MPCG method, it is crucial to select preconditioners with distinct properties, such that the preconditioned residuals do not become linearly dependent; otherwise, the matrix $\boldsymbol{P}_k^\top\boldsymbol{A}\boldsymbol{P}_k$ can become singular.
If the FE spaces in the multilevel hierarchy are $H_1$-orthogonal to each other, we can significantly enhance the search space, leading to improved convergence and also avoid the matrix $\boldsymbol{P}_k^\top\boldsymbol{A}\boldsymbol{P}_k$ from becoming singular.
In Algorithm~1, we can see the detailed algorithm for combining the additive MG and the MPCG method.

\begin{algorithm}[t]
	\label{algo:mpcg-addmg}
	\caption{An additive MG multipreconditioned CG algorithm}\label{alg:multi_cg}
	\SetKwData{Left}{left}\SetKwData{This}{this}\SetKwData{Up}{up}
	\SetKwComment{Comment}{$\triangleright$\ }{}
	\SetKwFunction{Union}{Union}\SetKwFunction{FindCompress}{FindCompress}
	\SetKwInOut{Input}{Input}\SetKwInOut{Output}{Output}
	\KwData{$\boldsymbol{A}, \boldsymbol{b},m, \nu, \boldsymbol{x}_0, (P_\ell^L)_{\ell=0,1,\ldots,L}, (R_L^\ell)_{\ell=0,1,\ldots,L}, (\pazocal{S}_\ell)_{\ell=0,1,\ldots,L}$}
	\KwResult{$\boldsymbol{x}_\ast = \boldsymbol{x}_k$}
	$\boldsymbol{r}_{0} \mapsfrom \boldsymbol{b} - \boldsymbol{A} \boldsymbol{x}_{0} $\Comment*[r]{Compute residual}
	$\boldsymbol{z}^\ell_{0} \mapsfrom P_\ell^L \pazocal{S}_\ell^\nu (R_L^\ell \boldsymbol{r}_{0}) \quad \forall \ell\in\{0,1,\ldots,L\}$ \Comment*[r]{Corrections from all levels}
	$\boldsymbol{P}_{0}=\boldsymbol{Z}_{0} \mapsfrom \begin{bmatrix}\boldsymbol{z}^0_{0}\ \vert\ \boldsymbol{z}^1_{0}\ \vert\ \cdots\ \vert\ \boldsymbol{z}^L_{0} \end{bmatrix}$ \Comment*[r]{Preconditioned residuals from $L$ levels}
	$k \mapsfrom 0 $\;
	\While{not converged}
	{
	$\boldsymbol{\alpha} \mapsfrom (\boldsymbol{P}_{k}^\top \boldsymbol{A} \boldsymbol{P}_{k})^{-1}(\boldsymbol{P}_{k}^\top \boldsymbol{r}_{k}) $ \Comment*[r]{Computing $\{\alpha_\ell\}_{\ell = 0, 1, \ldots,L}$ }
	$\boldsymbol{x}_{k+1} \mapsfrom \boldsymbol{x}_{k} + \boldsymbol{P}_{k} \boldsymbol{\alpha} $\Comment*[r]{updating the current iterate}
	$\boldsymbol{r}_{k+1} \mapsfrom \boldsymbol{r}_{k} - \boldsymbol{A} \boldsymbol{P}_{k} \boldsymbol{\alpha} $ \Comment*[r]{updating the residual}
	$\boldsymbol{z}^\ell_{k+1} \mapsfrom P_\ell^L \pazocal{S}_\ell^\nu (R_L^\ell \boldsymbol{r}_{k+1})  \quad \forall \ell\in\{0,1,\ldots,L\} $ \Comment*[r] {Corrections from all levels}

	$\boldsymbol{Z}_{k+1} \mapsfrom \begin{bmatrix}\boldsymbol{z}^0_{k+1}\ \vert \ \boldsymbol{z}^1_{k+1}\ \vert \ \cdots \ \vert \ \boldsymbol{z}^L_{k+1} \end{bmatrix}$ \Comment*[r]{Preconditioned residuals}
	$\boldsymbol{P}_{k+1} \mapsfrom \boldsymbol{Z}_{k+1} - \sum_{j=1}^{m} \boldsymbol{P}_j$ $ (\boldsymbol{P}_j^\top \boldsymbol{A} \boldsymbol{P}_j)^{-1}\boldsymbol{P}_j^\top \boldsymbol{A} \boldsymbol{Z}_{k+1} $ \Comment*[r]{Constructing $\beta$ }$k \mapsfrom k + 1$
	}
\end{algorithm}

\section{Numerical Experiments}
We investigate the performance of the additive MG-MPCG method using two numerical examples.

\begin{figure}[t]
	\centering
	\begin{subfigure}[b]{0.27\textwidth}
		\centering
    \includegraphics{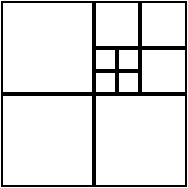}


		\caption{Sketch of fracture network}
		\label{fig:example2_sketch}
	\end{subfigure}
	\hfill
	\begin{subfigure}[b]{0.33\textwidth}
		\centering
		\includegraphics[scale=0.12]{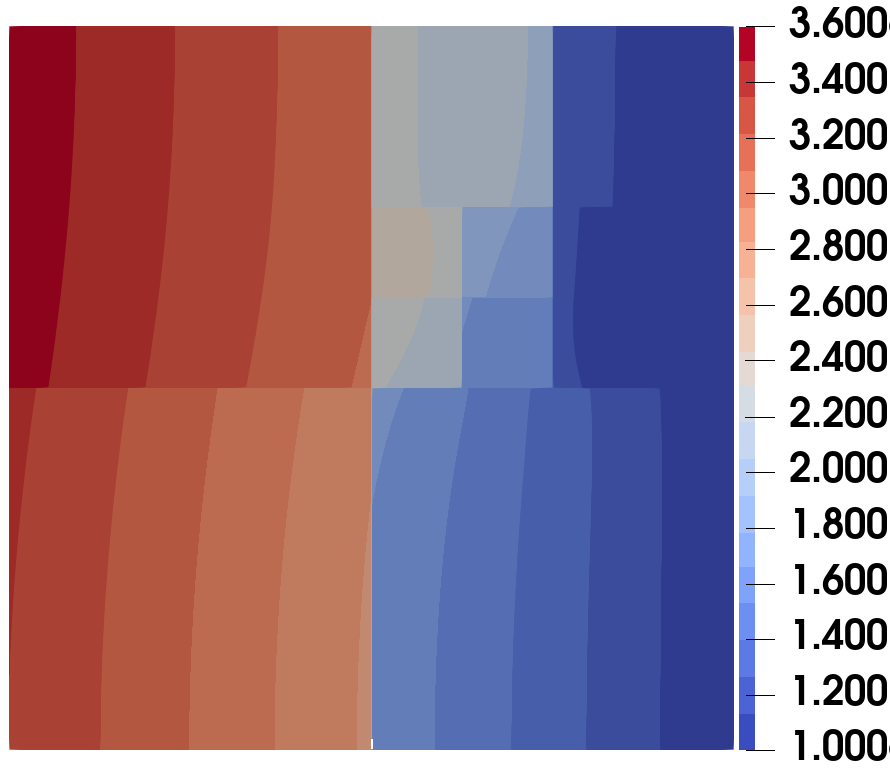}
		\caption{Blocking fractures:\\$k_m = 1$, $k_f = 10^{-4}$}
		\label{fig:example2_1e-4}
	\end{subfigure}
	\hfill
	\begin{subfigure}[b]{0.33\textwidth}
		\centering
		\includegraphics[scale=0.12]{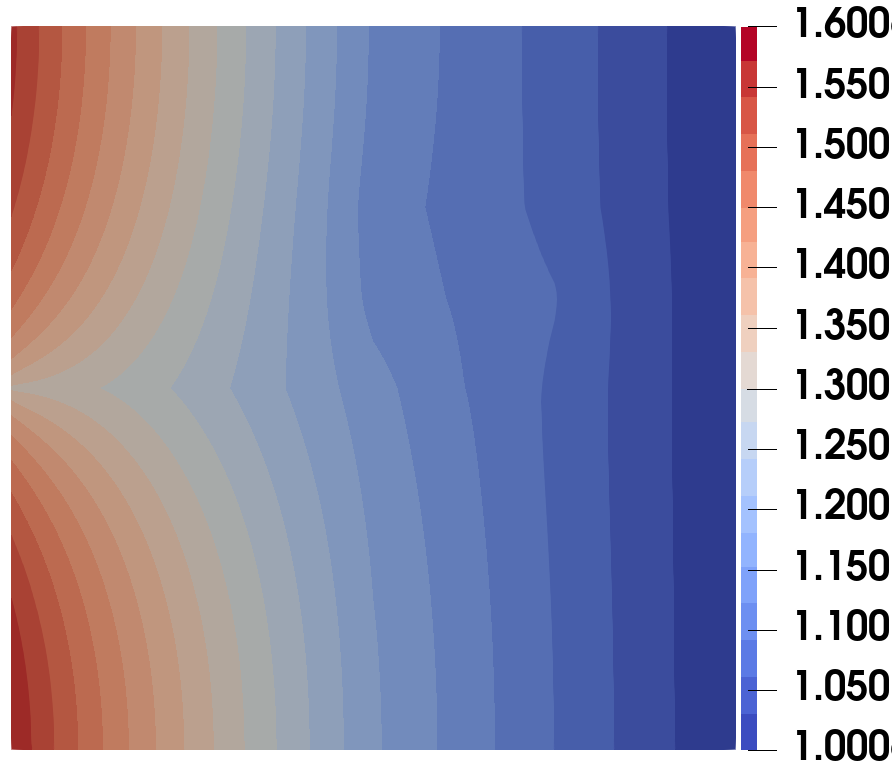}
		\caption{Conducting fractures:\\$k_m = 1$, $k_f = 10^{4}$}
		\label{fig:example2_1e4}
	\end{subfigure}
	\caption{Fracture network and pressure distribution for two cases (Example~2)}
	\label{fig:example2}
\end{figure}

\vspace{0.1cm}
\noindent\emph{Example 1: Anisotropic diffusion:}\label{para:example1}
We consider a domain $\Omega = (0,1)^2$ with anisotropic diffusion.
The values considered for $k_{xx}$ range across several orders of magnitude, $k_{xx} \in \{10^n\}_{n=-7,-6,\ldots,0}$, while we set $k_{yy} = 1$.
We utilize a mesh with $160 \times 160$ elements, which we set as the finest level in the multilevel hierarchy.
In the multigrid method, we consider a nested hierarchy of meshes with $4$ levels, while on the coarsest level, we have only $20 \times 20$ elements.

\vspace{0.1cm}
\noindent\emph{Example 2: Pressure distribution in the fracture network:}
In this example, we consider the domain $\Omega = (0,1)^2$ to model pressure distribution in a fracture network, considering various flow dynamics parameters~\cite{11_nestola_novel_2020}.
The fracture network is characterized by a uniform fracture thickness ($\delta$) of $10^{-4}$.
The matrix permeability ($k_m$), a measure of the ability of the matrix to transmit fluids, is set to $1$.
The fracture permeability ($k_f$) is chosen as $k_f \in \{10^n\}_{n=-4,-3,\ldots,4}$, ranging from the blocking to conducting phase.
The model involves approximately $3 \cdot 10^{5}$ degrees of freedom (DoFs).
We enforce $g_N = 0$ as the Neumann boundary condition on the top and bottom edge of the domain, while $g_N = 1$ on the left part of the boundary and $u_D = 1$ is prescribed on the right edge of the domain.
A sketch of the fracture network can be seen in Figure~\ref{fig:example2_sketch}, while the fracture flow in the network for the blocking fractures can be seen in Figure~\ref{fig:example2_1e-4} and conducting fractures in Figure~\ref{fig:example2_1e4}.
An Algebraic Multigrid (AMG) method using BoomerAMG from the Hypre framework~\cite{8_henson_boomeramg_2002} constructs the multilevel hierarchy.
The solution methods follow the strategies detailed in Example~1.

\vspace{0.1cm}
\noindent\emph{Setup for the solution methods:}
We use the symmetric Successive Over-Relaxation (SSOR) method as a smoother with $3$ pre- and $3$ post-smoothing steps for the multiplicative MG while $6$ smoothing steps on each level for the additive MG, ensuring the same amount of work is done at each level for both multigrid variants.
Furthermore, for the MPCG method, we store the last $5$ pairs of vectors to estimate the value of the parameter $\beta$, balancing accuracy and memory use for computational efficiency.
While storing more vectors can improve convergence by preserving the conjugacy of the search directions, it becomes memory-intensive, requiring storage proportional to {${(L+1)\times n \times m}$}, increasing $m$ enhances convergence but raises memory costs.
We store the last $5$ pairs of vectors, as it was shown in\cite{4_bridson_multipreconditioned_2006} that a few vectors are sufficient and increasing the number of vectors further might not significantly improve convergence.

\subsection{Comparison with the other solution methods} %

\begin{figure}[t]
\includegraphics{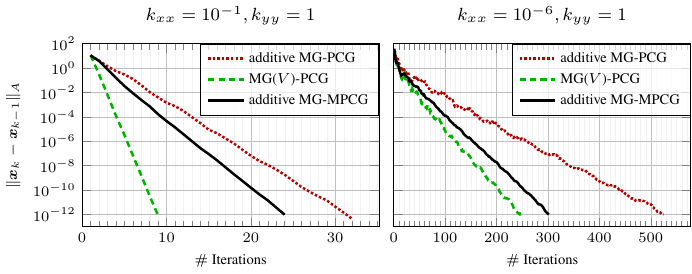}
	\caption{Iterations $v/s$ corrections for two different test cases of Example$~1$}
	\label{fig:example1_convergence}
\end{figure}

In this section, we compare the convergence behavior of the additive MG-MPCG method, the additive MG-PCG method, and the multiplicative MG PCG method.

For the Example~1, we consider a high anisotropic case ($k_{xx}=10^{-6}$, $k_{yy}=1$) and a low anisotropic case ($k_{xx}=10^{-1}$, $k_{yy}=1$).
As shown in Figure~\ref{fig:example1_convergence}, in both cases, the MPCG method outperforms the standard PCG method when preconditioned with the additive MG method.
However, as expected, the multiplicative MG-PCG method outperforms both methods.
For Example~1, we have employed a multigrid method with geometric hierarchy using the uniform refinement strategy.
Hence, we can see that the multilevel hierarchy is not optimal for the highly anisotropic case.
Consequently, we observe that the number of iterations required for convergence increases when anisotropy is considered.
For the highly anisotropic case, the performance of the additive MG-MPCG method is closer to the MG($V$)-PCG method.

To investigate this behavior further, we carried out the same experiment with more cases with slowly increasing anisotropy.
Figure~\ref{fig:compare_coeffs} (left) shows that problems with reduced anisotropy converge significantly faster than those with high anisotropy.
Among the tested methods, the additive MG-PCG requires the highest number of iterations for convergence, and this requirement amplifies more noticeably with increasing anisotropy.
Interestingly, the additive MG-MPCG method, while exhibiting a similar trend, does not experience a significant increase in the number of iterations
Its convergence rate grows at a rate comparable to the MG($V$)-PCG methods.

\begin{figure}[t]
\includegraphics{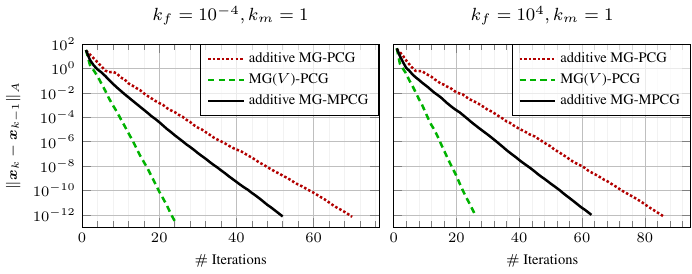}

	\caption{Iterations $v/s$ corrections for two different test cases of Example$~2$}
	\label{fig:example2_convergence}
\end{figure}

\begin{figure}[t]
\includegraphics{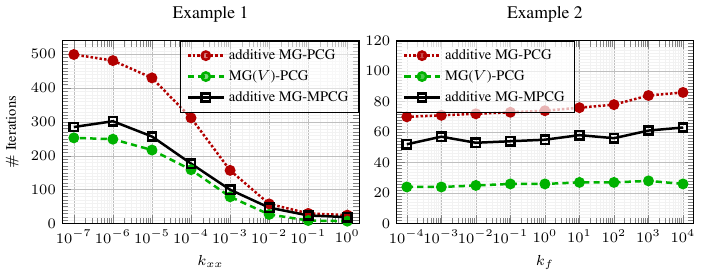}

	\caption{Comparing the number of iterations required for convergence for various values of the diffusion coefficients $k_{xx}$ and fracture permeability $k_f$.}
	\label{fig:compare_coeffs}
\end{figure}

For Example~2, we consider two different test cases with a varying range of fracture permeabilities, where the two extreme cases are chosen, the blocking fractures ($k_m = 1, k_f = 10^{-4}$) and the conducting fractures ($k_m = 1, k_f = 10^4$).
In Figure~\ref{fig:example2_convergence}, we observe the convergence of the three solution methods for the blocking and conductive fractures.
Similar to the previous example, the additive MG-MPCG method outperforms the additive MG-PCG method for both test cases.
In this example, we utilize BoomerAMG to build the multilevel hierarchy, which takes into account the permeabilities of the fractures when constructing the transfer operators.
Such agglomeration based construction enhances the quality of the multilevel hierarchy compared to that in Example 1, resulting in a more effective hierarchy of the subspaces.
We can observe the effect of the transfer operators in Figure~\ref{fig:compare_coeffs} (right).
In this case, the number of iterations for additive MG-PCG and additive MG-MPCG methods remains relatively constant across different permeabilities, with the additive MG-MPCG method showing a slight increase in iterations as permeability increases.

Optimal transfer operators and smoothers typically ensure a constant number of iterations for additive or multiplicative  multilevel methods regardless of problem size.
However, our anisotropic problem with non-optimal transfer operators results in increased iterations as problem size grows. Our experiments indicate that the trend observed in Figure 4 will persist with increasing problem size.

\subsection{Comparing components of $\boldsymbol{\alpha}$} %
\begin{figure}[t]
	\begin{subfigure}[b]{0.49\textwidth}
		\centering
		\includegraphics[scale=0.193]{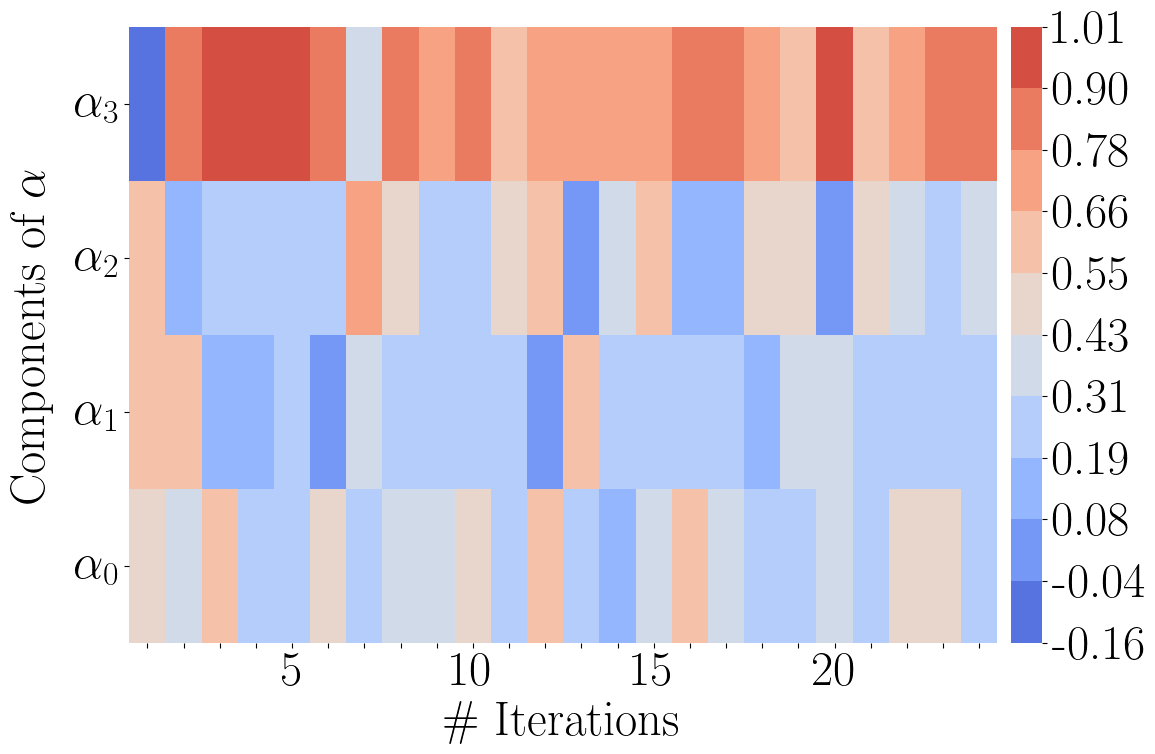}
		\caption{$k_{xx} = 10^{-1}$, $k_{yy} = 1$}
		\label{fig:alphas_example1_1e-1}
	\end{subfigure}
	\hfill
	\begin{subfigure}[b]{0.49\textwidth}
		\centering
		\includegraphics[scale=0.193]{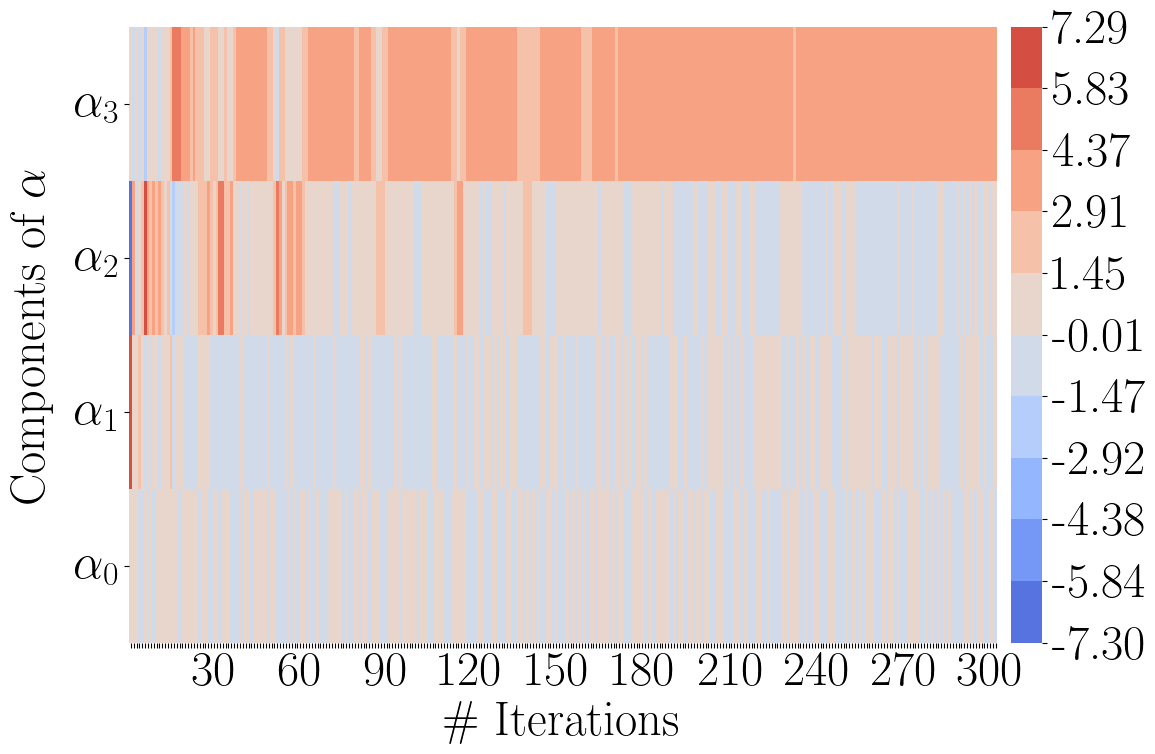}
		\caption{$k_{xx} = 10^{-6}$, $k_{yy} = 1$}
		\label{fig:alphas_example1_1e-6}
	\end{subfigure}
	\caption{Example~1: Heatmap of components of $\alpha_\ell$ associated with different level of the multilevel hierarchy}
	\label{fig:alphas_example1}
\end{figure}

\begin{figure}[t]
	\begin{subfigure}[b]{0.49\textwidth}
		\centering
		\includegraphics[scale=0.192]{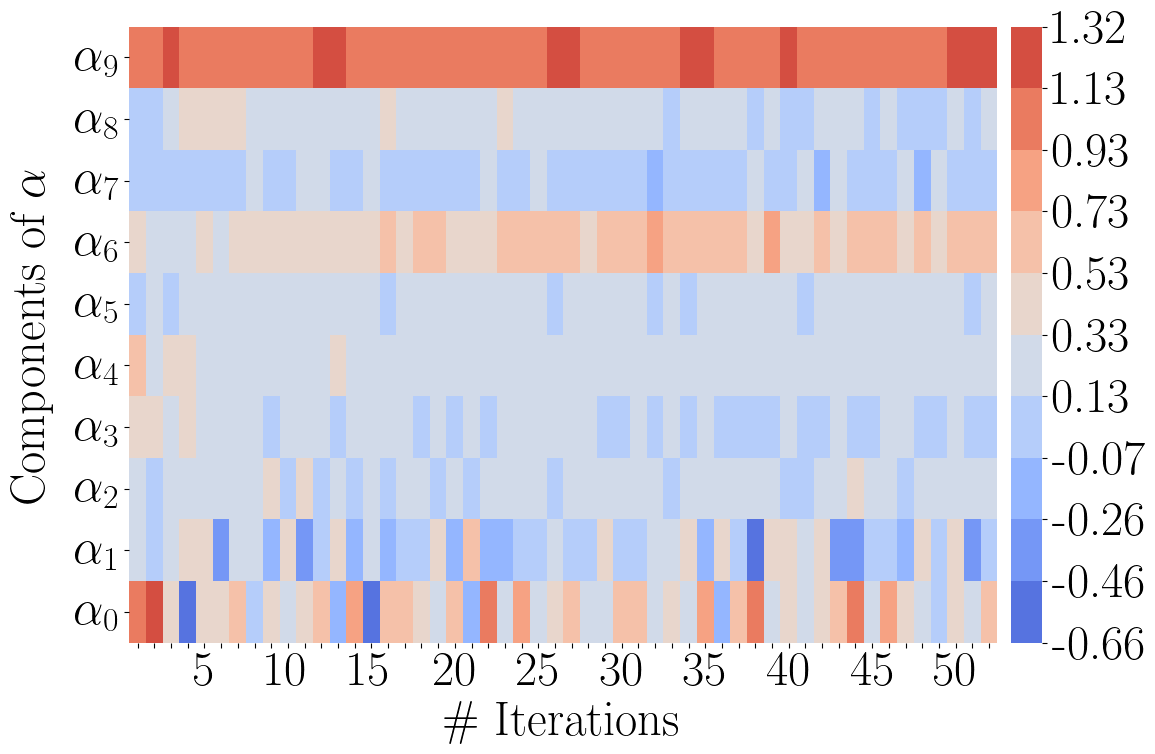}
		\caption{$k_{f} = 10^{-4}$, $k_{m} = 1$}
		\label{fig:alphas_example2_1e-4}
	\end{subfigure}
	\hfill
	\begin{subfigure}[b]{0.49\textwidth}
		\centering
		\includegraphics[scale=0.192]{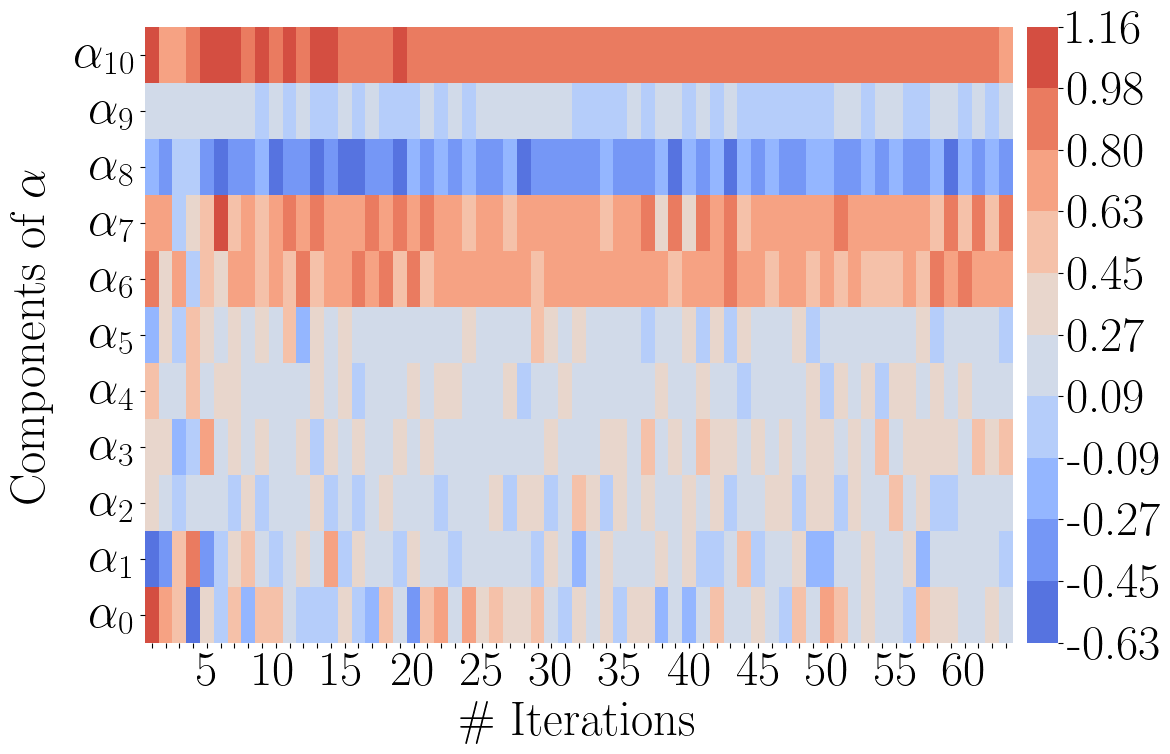}
		\caption{$k_{f} = 10^{4}$, $k_{m} = 1$}
		\label{fig:alphas_example2_1e4}
	\end{subfigure}
	\caption{Example~2: Heatmap of components of $\alpha_\ell$ associated with different level of the multilevel hierarchy}
	\label{fig:alphas_example2}
\end{figure}

In this section, we study the components of $\boldsymbol{\alpha}$ that are used to construct the optimal search direction in the MPCG method over all iterations.
Figure~\ref{fig:alphas_example1} shows the heatmap of values of $\boldsymbol{\alpha}$ associated with each level of the multilevel hierarchy for the two test cases of Example~1, while Figure~\ref{fig:alphas_example2} demonstrates the heatmap of values of $\boldsymbol{\alpha}$ for the blocking and conducting fracture cases of Example~2.
From both figures, it is interesting to note that the values of $\alpha_\ell$ can also become negative.

From Figure~\ref{fig:alphas_example1}, we can see that for the few initial iterations for both the low anisotropic and the high anisotropic cases, the values of $\alpha_\ell$ are positive.
However, after a few iterations, the values of $\alpha_\ell$ associated with the coarser levels become smaller.
This behavior is prominent in the high anisotropic case (c.f.~Figure~\ref{fig:alphas_example1_1e-6}) as only positive values of $\alpha_\ell$ are associated with the finest level in the latter iterations.
This suggests that the corrections from the coarser level are not optimal, and the MPCG method automatically discards the correction from the coarse level by attributing smaller values of $\alpha_\ell$.

In Example 2, we observe similar behavior in the values of $\alpha_\ell$.
Notably, some levels consistently exhibit negative values of $\alpha_\ell$.
Despite this, the finest level consistently shows the largest values of $\alpha_\ell$ in both the blocking and conducting cases, and these values are significantly higher compared to other levels.
This observation indicates that some levels in the hierarchy might not contribute positively, as suggested by the values of $\boldsymbol{\alpha}$.

Overall, the MPCG method demonstrates a notable capability in filtering and prioritizing levels within the hierarchy.
It effectively distinguishes the most beneficial corrections for the convergence process, ensuring that even levels appearing less useful at a glance are integrated into the calculation to achieve the desired smoothing across the spectrum.

\section{Conclusion}
In conclusion, we present a novel approach for integrating the additive multigrid method with  the MPCG method.
The numerical studies lead to the proposition that, especially in situations where an optimal hierarchy of multilevel subspaces is unavailable, the additive MG-MPCG method can be highly beneficial.
Its robust convergence properties and additive nature indicate a strong potential for outperforming the multiplicative MG-PCG method, notably in parallel computing environments.
Future work will focus on leveraging parallel implementations, and combining the proposed methodology with a matrix-free strategy~\cite{9_kothari_multigrid_2022}.
Furthermore, this work can be seamlessly extended within the unfitted finite element framework for interface problems in the domain decomposition framework~\cite{10_kothari_multigrid_2021*a}.

\begin{acknowledgement}
	This work is supported by the Swiss National Science Foundation (SNSF) and the Deutsche Forschungsgemeinschaft for their through the project SPP 1962 ``Stress-Based Methods for Variational Inequalities in Solid Mechanics:~Finite Element Discretization and Solution by Hierarchical Optimization" [186407].
	We also acknowledge the support of Platform for Advanced Scientific Computing (PASC) through the project FraNetG:~Fracture Network Growth.
\end{acknowledgement}

\end{document}